\documentclass[12pt]{amsart}
\usepackage{amscd,amssymb,amsmath}
\newtheorem{thm}{Theorem}[section]

\theoremstyle{definition}

\theoremstyle{remark}

\numberwithin{equation}{section}

\def\Is{{\mathrm{Is}}\,}

\def\QED{\nobreak\quad\ifmmode\roman{Q.E.D.}\else{\rm Q.E.D.}\fi}

\def\R{{\mathbf R}}

\def\G{\Gamma}

\def\sB{{\mathcal B}}

\def\sF{{\mathcal F}}

\def\sU{{\mathcal U}}

\def\sbs{\subset}
\def\rar{\rightarrow}

\def\ti{\times}
\def\obr{^{-1}}
\def\stm{\setminus}

\begin{document}

\title[The Urysohn space is homeomorphic to a Hilbert space]
{The Urysohn universal metric space is homeomorphic to 
a Hilbert space}

\author[V. Uspenskij]{Vladimir Uspenskij}

\address{321 Morton Hall, Department of Mathematics, Ohio
University, Athens, Ohio 45701, USA}

\email{uspensk@math.ohiou.edu}


\date{3 September 2003}


\begin{abstract} 
The Urysohn universal metric space $U$ is characterized 
up to isometry by the following properties: (1) $U$ is 
complete and separable; (2) $U$ contains an isometric copy
of every separable metric space; (3) every isometry between
two finite subsets of $U$ can be extended to an isometry
of $U$ onto itself. We show that $U$ is homeomorphic to
the Hilbert space $l_2$ (or to the countable power of the
real line). 
\end{abstract}

\maketitle

\setcounter{tocdepth}{1}

\section{Introduction} 
The Urysohn universal metric space $U$ is characterized 
up to isometry by the following properties: (1) $U$ is 
complete and separable; (2) $U$ contains an isometric copy
of every separable metric space; (3) every isometry between
two finite subsets of $U$ can be extended to an isometry
of $U$ onto itself. (An {\it isometry} is a distance-preserving
bijection; an {\it isometric embedding} is a distance-preserving
injection.) The aim of the present paper is to show
that the Urysohn space $U$ is homeomorphic to a Hilbert space
(equivalently, to the countable power of the real line).
This answers a question raised by Bogaty\u\i, Pestov and Vershik.

There is another characterization of $U$.
Let us say that a metric space $M$ is 
{\it injective with respect to finite spaces}, or 
{\it finitely injective} for short, if for every finite metric space
$L$, every subspace $K\sbs L$ and every isometric embedding $f:K\to M$
there exists an isometric embedding $\bar f:L\to M$ which extends $f$.
Define {\it compactly injective} metric spaces similarly,
considering compact (rather than finite) spaces $K$ and $L$.
If a metric space $M$ contains an isometric copy of every finite
metric space and satisfies the condition (3) above, then $M$ 
is finitely injective. 
Indeed, given finite metric spaces
$K\sbs L$ and an isometric embedding $f:K\to M$, we find an isometric
embedding $g:L\to M$ and extend the isometry $gf\obr:f(K)\to g(K)$
to an isometry $h$ of $M$ onto itself. Then $h\obr g:L\to M$
is an isometric embedding of $L$ which extends $f$.
Conversely, let $M$ be a finitely injective
metric space. Then every countable metric space admits an isometric
embedding into $M$ (use induction). If $M$ is also complete, it follows
that $M$ contains an isometric copy of every separable metric space.
Assume now that $M$ is also separable, and let $f:K\to L$ be an isometry
between two finite subsets of $M$. Enumerating points of a dense
countable subset of $M$ and using the back-and-forth method we can extend
$f$ to an isometry between two dense countable subsets of $M$ and
then to an isometry of $M$ onto itself. The same argument shows
that any two complete separable finitely injective metric spaces
are isometric. Thus the Urysohn space $U$ 
is the unique (up to isometry) metric space
which is complete, separable and finitely injective.

The existence of $U$ was proved by
Urysohn \cite{Ur25, Ur27}. An easier construction was found some 50 years
later by Kat\v etov \cite{Kat}, who also gave an example of a non-complete
separable metric space satisfying the conditions (2) and (3)
above, thus answering a question of Urysohn. Kat\v etov's construction
was used in \cite{Usp1, Usp2, Usp3} to prove that the topological group
$\Is(U)$ of all isometries of $U$ is universal, in the sense that
it contains an isomorphic copy of every topological group
with a countable base. A deep result concerning the group $G=\Is(U)$
was established by V.Pestov: the group $G$ is extremely amenable,
i.e., every compact space with a continuous action of $G$ has 
a $G$-fixed point \cite{P1, KPT}.

A.M.Vershik showed that the space $U$ can be obtained
as the completion of a countable metric space equipped with a 
metric which is either``random"  or generic in the sense of Baire
category \cite{Versh1, Versh2, Versh3}. 

Bogaty\u\i{} \cite{B} 
proved that any isometry between two compact subsets of $U$
can be extended to an isomerty of $U$ onto itself. 
It follows by the same argument that we used above for finitely
injective spaces
that $U$ is compactly injective (and is the unique complete separable
compactly injective metric space). 
Using this, we deduce our 
Main Theorem from Toru\'nczyk's Criterion \cite{Torun, Revis}:
a complete separable metric space
$M$ is homeomorphic to the Hilbert space $l_2$ 
if and only if $M$ is AR (= absolute
retract) and has the discrete approximation property
(this notion is defined below). 
Recall that all infinite-dimensional separable Banach spaces
are homeomorphic to each other and to the countable power of the
real line.

Given an open cover $\sU$ of a space $X$, two points $x,y\in X$
are said to be {\it $\sU$-close} if there exists $U\in \sU$
such that $x,y\in U$. 
A family of subsets of a space
$X$ is {\it discrete} if every point in $X$ has a neighbourhood which
meets at most one member of the family.
A metric space $M$ has the {\it discrete approximation property}
if for every sequence $K_1, K_2, \dots$ of compact subspaces
of $M$ and every open cover $\sU$ of $M$ there
exists a sequence of maps $f_i:K_i\to M$ such that for every $i$
and every $x\in K_i$ the points $x$ and $f_i(x)$ are $\sU$-close
and the sequence 
$(f_i(K_i))$ is discrete. 
Equivalently \cite{Torun}, a metric space 
$(M,d)$ has the discrete approximation property
if and only if 
for every sequence $K_1, K_2, \dots$ of compact subspaces      
of $M$ 
and every continuous function $h$ on $M$ with values $>0$
there exists a sequence of maps $f_i:K_i\to M$ such that 
$d(x,f_i(x))\le h(x)$ for every $i$ and every $x\in K_i$,
and the sequence 
$(f_i(K_i))$ is discrete. 

Let us reformulate Toru\'nczyk's Criterion in the form that is convenient
for our purposes.
We say that a topological space $X$ is {\it homotopically trivial}
if $X$ has trivial homotopy groups, that is, every map of the $n$-sphere
$S^n=\partial B^{n+1}$ to $X$ admits an extension over the $(n+1)$-ball
$B^{n+1}$ ($n=0,1,\dots$). (The term {\it weakly homotopically trivial}
might be more appropriate.) Every contractible space is
homotopically trivial; the converse 
in general is not true. The empty space is homotopically trivial.
If a metric space $M$
has a base $\sB$ such that for every non-empty finite subfamily
$\sF\sbs\sB$ the intersection $\cap\sF$ is homotopically trivial,
then $M$ is ANR \cite[Theorem 5.2.12]{vM}. 
A metric space is AR if and only if
it is homotopically trivial and ANR \cite[Theorem 5.2.15]{vM}.
Thus Toru\'nczyk's
Criterion can be reformulated as follows:

\begin{thm}[{\bf Toru\'nczyk's Criterion}] A complete separable metric
space $M$ is homeomorphic to a Hilbert space if and only if the 
following conditions hold:
\begin{enumerate}
\item[(i)] there is a base $\sB$ for $M$ such that $U,V\in\sB$ implies
$U\cap V\in \sB$, and every $U\in \sB$ is homotopically trivial;
\item[(ii)] $M$ is homotopically trivial;
\item[(iii)] $M$ has the discrete approximation property.
\end{enumerate}
\end{thm}

In the next section we show that the space $U$ satisfies the conditions
of this criterion.
 
\section{Proof of the main theorem}

\begin{thm}[Main Theorem]
The Urysohn universal space $U$ is homeomorphic to a Hilbert space.
\end{thm}

\begin{proof}
We check the three conditions of Toru\'nczyk's criterion.

(a) 
Let $\sB$ be the base for $U$ consisting of all open balls $O(a,r)=\{x\in U:
d(x,a)<r\}$ and their finite intersections. We claim that every member
$V=\cap_{i=1}^k O(a_i, r_i)$ of $\sB$ is homotopically trivial.
Let a map $f:S^n\to V$ be given. We must construct an extension
$\bar f:B^{n+1}\to V$.

Every metric space admits an isometric embedding into a normed linear
space. Thus we may consider $U$ as a subspace of a Banach space $B$.
Let $V'=\cap_{i=1}^k O'(a_i, r_i)$, where $O'(a,r)$ is the open ball
centered at $a$ of radius $r$ in the space $B$. Then $V=V'\cap U$.
Being a convex subset of a normed linear space, the space $V'$ is 
contractible (in fact it is AR \cite[Theorem 1.4.13]{vM}), 
so the map $f:S^n\to V$ can be extended to a map 
$g:B^{n+1}\to V'$. Let $A=\{a_1,\dots, a_k\}$, $K=f(S^n)\cup A$
and $L=g(B^{n+1})\cup A$. Then $K$ and $L$ are compact, $K\sbs L\cap U$.
Since $U$ is compactly injective, the identity map of $K$ can be extended
to an isometric embedding $h:L\to U$. Let $\bar f=hg:B^{n+1}\to U$.
Then $\bar f$ extends $f$. The range of $\bar f$ is contained in $V$,
since for every $x\in B^{n+1}$ and $i=1,\dots, k$ we have 
$d(\bar f(x),a_i)=d(h(g(x)),h(a_i))=d(g(x),a_i)<r_i$
(note that $h(a_i)=a_i$, since $a_i\in K$ and $h$ fixes all points in $K$).

(b) The space $U$ is homotopically trivial. The proof is the same
as above but easier, since we do not have to care about points
$a_1,\dots, a_k$.

(c) We prove that $U$ has the discrete approximation property.
Let $K_1,\dots, K_n, \dots$ be a sequence of non-empty 
compact subsets of $U$,
and let $h$ be a continuous function on $U$ with values $>0$.
We must construct a discrete sequence $(L_n)$ of compact subsets
of $U$ and a sequence of maps $f_n:K_n\to L_n$ such that 
$d(f_n(x),x)\le h(x)$ for every $n\ge1$ and $x\in K_n$.

We'll need the notion of union of two metric
spaces with a subspace amalgamated. Suppose that $M_1, M_2, A$
are metric spaces, $A\ne\emptyset$, and isometric embeddings
$f_i:A\to M_i$, $i=1,2$, are given. The union $M$ of $M_1$ and $M_2$
with the subspace $A$ amalgamated is characterized by the following
properties: there exist isometric embeddings $h_i:M_i\to M$
such that $M=h_1(M_1)\cup h_2(M_2)$, $h_1f_1=h_2f_2$, and for every
$x\in M_1\stm f_1(A)$, $y\in M_2\stm f_2(A)$ 
$$
d(h_1(x),h_2(y))=\inf\{d_1(x,f_1(z))+d_2(f_2(z),y):z\in A\},
$$
where $d,\ d_1,\ d_2$ are the metrics on $M,\ M_1,\ M_2$, respectively.
It is easy to see that such a space $M$ exists and in the obvious sense
is unique.

Let $N_i\sbs K_i\ti \R$ be the union of $K_i\ti\{0\}$ and 
the graph of the restriction of $h$ on $K_i$. Equip $K_i\ti\R$
with the metric $\rho$ defined by 
$$
\rho((x,t),(y,s))=d(x,y)+|s-t|,
$$
and consider the induced metric on $N_i$.

We now construct a sequence $(L_n)$ of compact subsets of $U$ by
induction. Suppose the sets $L_i$ have been defined for $i<n$.
We define $L_n$. Consider two compact metric spaces: 
$K_n\cup\bigcup_{i<n}L_i$ and $N_n$. Since $K_n$ lies in the first
space and has a natural embedding into the second one
(we mean the embedding $x\mapsto (x,0)$), 
we can construct their union with the subspace
$K_n$ amalgamated. Write this union as 
$P=\bigcup_{i<n}L_i\cup K_n\cup\G_n$, where
$\G_n=\{(x,h(x)):x\in K_n\}$ is the graph of $h\restriction K_n$.
Since $U$ is compactly injective, there exists an isometric
embedding $\phi:P\to U$ which is identity on each $L_i$ ($i<n$)
and on $K_n$. Let $L_n=\phi(\G_n)$. Let $f_n:K_n\to L_n$ be
the composition of the map $x\mapsto(x,h(x))$ 
from $K_n$ onto $\G_n$ and $\phi$.
For every $x\in K_n$ 
the distance from $(x,0)$ to $(x,h(x))$ in $N_n$ is equal to
$h(x)$, hence the distance from $x$ to $(x,h(x))$ in $P$ and
the distance from $x$ to $f_n(x)$ in $U$ also are equal to
$h(x)$. Thus $f_n$ moves every $x\in K_n$ by $h(x)$.

Note that for every $x\in K_n$ and $y\in \bigcup_{i<n}L_i$
the distance from $f_n(x)$ to $y$ is $\ge h(x)$. 
Indeed, by our construction
this distance is equal to the distance from $(x,h(x))\in \G_n$
to $y$ in $P$ and thus also to 
\begin{align*}
\inf\{d(y,z)+\rho((z,0),(x,h(x))):z\in K_n\}&=
\inf\{d(y,z)+d(z,x)+h(x):z\in K_n\}\\
&=d(y,x)+h(x)\ge h(x).
\end{align*}

To conclude the proof, we must show that the sequence $(L_n)$
is discrete. Assume the contrary. Since the sequence $(L_n)$
is disjoint, there exists an infinite set $A$ of positive integers
and points $y_i\in L_i$ ($i\in A$) such that the sequence
$\{y_i:i\in A\}$ converges to some $p\in U$. Write $y_i=f_i(x_i)$,
where $x_i\in K_i$. The distance from $y_n$ to $\{y_i:i<n,\ i\in A\}
\sbs \bigcup_{i<n}L_i$ tends to zero as $n\in A$ tends to infinity. 
On the other hand, we saw in the preceding paragraph that this distance
is $\ge h(x_n)$. Therefore the sequence $\{h(x_n):n\in A\}$ tends
to zero. Since $d(x_n, y_n)=d(x_n, f_n(x_n))=h(x_n)\rar 0$
and $y_n\rar p$, it follows that $x_n\rar p$.
But this contradicts the continuity of $h$ at $p$: we have $h(p)>0$,
$x_n\rar p$ and $h(x_n)\rar 0$.
\end{proof}

\end{document}